\documentclass[preprint,12pt,number]{elsarticle}
\usepackage{xcolor,amsmath,amssymb,amsfonts,array,graphicx}
\usepackage{hyperref,fontenc,lmodern}
\hypersetup{colorlinks=false,allbordercolors=blue,pdfborderstyle={/S/U/W 1}}
\usepackage[utf8]{inputenc}
\usepackage[T1]{fontenc}
\usepackage[polish,english]{babel}
\usepackage{hyperref,float}
\AtBeginDocument{\selectlanguage{english}}

\journal{Fractals}

\begin{document}
	
	\begin{frontmatter}
		
		\title{Fractal Patterns in Discrete Laplacians: Iterative Construction on 2D Square Lattices}
		
		\author{Ma\l gorzata Nowak-K\c epczyk} 
		
		\affiliation{organization={Faculty of Natural and Technical Sciences, The John Paul II Catholic University of Lublin},
			addressline={Konstantynow 1H}, 
			city={Lublin},
			postcode={20-708}, 
			country={Poland, email: malnow@kul.pl}}
		
\begin{abstract}
We investigate the iterative construction of discrete Laplacians on 2D square lattices, revealing emergent fractal-like patterns shaped by modular arithmetic. While classical 2222-style iterations reproduce known structures such as the Sierpiński triangle, our alternating binary–ternary (2322-style) process produces a novel class of aperiodic figures. These display low density variance, minimal connectivity loss, and non-repetitive organization reminiscent of Dekking’s sequences. Fourier and autocorrelation analyses confirm their quasi-periodic nature, suggesting structural analogies with self-assembly processes and fractal-based design.

 The results highlight the structural richness of modular Laplacian dynamics and motivate further study of non-periodic discrete systems.

\end{abstract}

\vspace{6mm}
		
		\begin{graphicalabstract}\centering
			\includegraphics[width=0.9\textwidth]{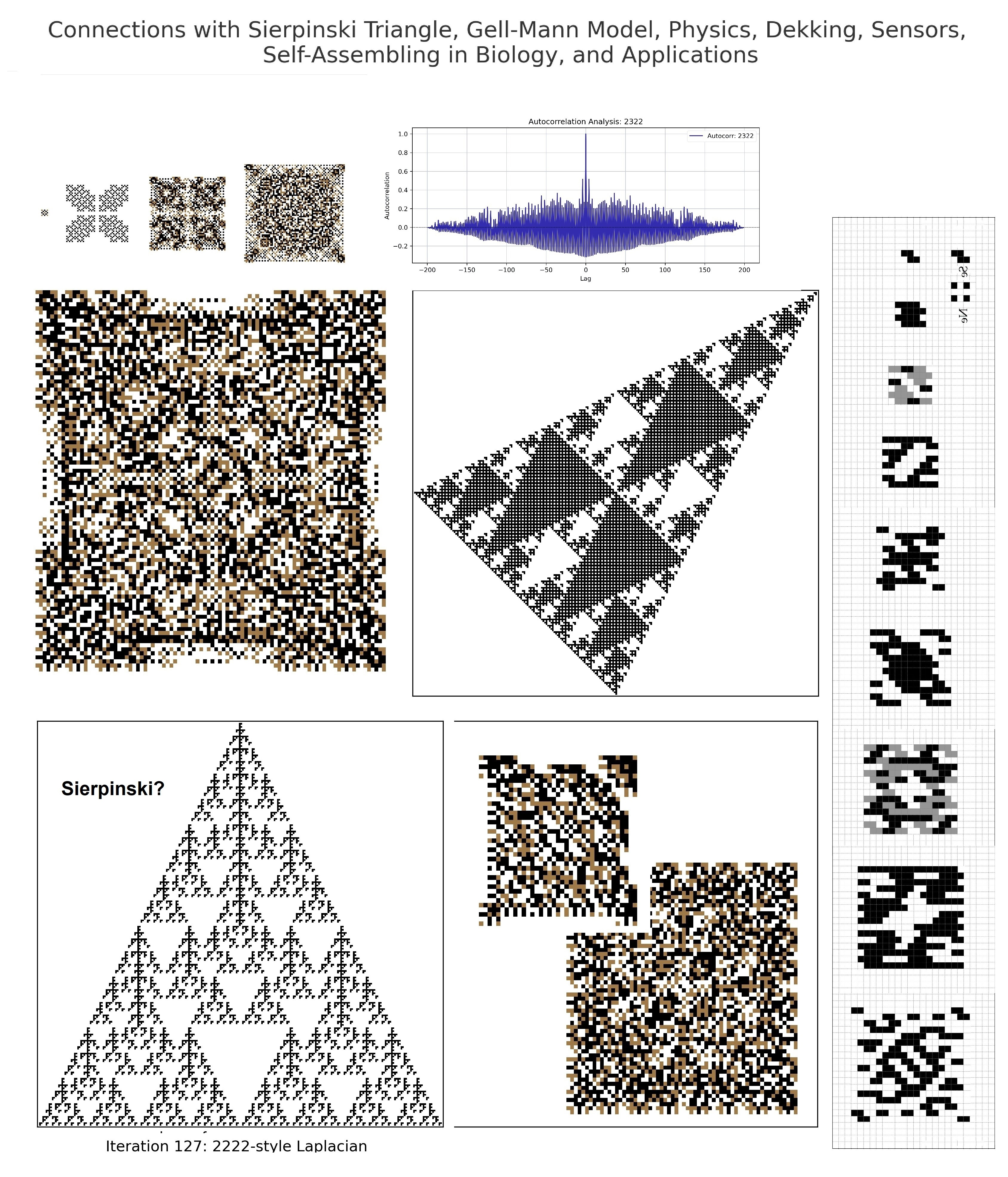}
		\end{graphicalabstract}

\begin{keyword}discrete Laplacian, modular arithmetic, fractals, aperiodic patterns, 2D lattices, Sierpiński triangle, Dekking sequence, sensor networks
\end{keyword}

\end{frontmatter}


\vspace*{0.05in}

\section*{Introduction}
The discrete Laplacian is a fundamental mathematical operator with applications in physics, biology, and computational science. It plays a crucial role in graph theory, numerical analysis, and dynamical systems, serving as a discrete analog to the continuous Laplace operator. Iterating the discrete Laplacian on 2D square lattices gives rise to complex emergent structures, often displaying fractal-like organization.

Previous studies of Discrete Laplacians on 2D lattices have documented visually striking fractal patterns, such as snowflake-like and carpet-like motifs, arising from binary iterations alone \cite{Aiba}. These patterns emerge due to modular arithmetic constraints and exhibit diverse forms, depending on the chosen seed and neighborhood structure. However, binary (modulo 2) iterations tend to produce periodic structures and, at specific iterations, dissociate into their initial seed configurations. This motivates the exploration of higher-order modular arithmetic, such as ternary (modulo 3) and quaternary (modulo 4) iterations, which significantly expand the space of possible patterns and often result in aperiodic, structurally richer formations.

Earlier investigations of modular Laplacian dynamics 
in the purely binary setting are presented in \cite{NowakKepczyk2015}.

In this paper, we investigate a special class of iterative Laplacian constructions, focusing on the interplay between binary and ternary arithmetic in pattern formation. We introduce the 2322-style pattern, an alternating sequence of binary and ternary iterations, which gives rise to structurally distinct, non-repetitive figures. Unlike purely binary 2222-style figures, which exhibit dissociative periodicity, 2322-style figures display low density variance, minimal connectivity loss, and quasi-aperiodicity. These properties are reminiscent of Dekking’s non-repetitive sequences, which play a significant role in combinatorial mathematics and aperiodic order.

Beyond their mathematical significance, these structures may have broader structural implications. The concluding section of this paper will discuss their potential relevance, including connections to biological self-organization, sensor technology, material science, and fractal-based computational algorithms.

The paper is structured as follows: Section 1 introduces the iterative construction process; Section 2 explores binary figures and Sierpiński-like patterns; Section 3 presents the distinguishing properties of 2322-style figures; Section 4 compares cell sequence statistics, including Fourier and autocorrelation analysis; and Section 5 discusses applications and outlines directions for future research.


\section{Iterative Dynamical Systems of Discrete Laplacians}

\noindent

\subsection{Initial Conditions and Seed Configurations} 
We consider an automaton initialized with binary values:
\[
\nu_0(p) =
\begin{cases} 
	1, \quad &\text{if } p \text{ belongs to the seed,} \\
	0, \quad &\text{otherwise}.
\end{cases}
\]
where 1 represents an occupied cell ($\rule{7pt}{7pt}$) and 0 denotes an unoccupied cell ($\square$). 

The choice of seed configuration significantly influences the evolutionary dynamics of the system. Some commonly used seeds include:
\begin{itemize}
	\item Single point (minimal seed, simplest growth),
	\item Line segment (introducing directional spread),
	\item Geometric clusters (leading to more intricate self-organization).
\end{itemize}

\begin{figure}[h!]
	\centering
	\includegraphics[width=0.8\textwidth]{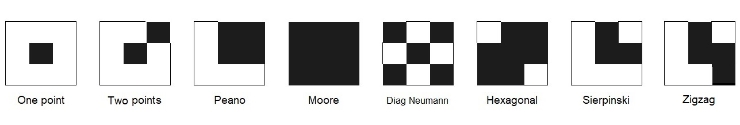}
	\caption{Examples of seed configurations.}
	\label{fig:seeds}
\end{figure}

\subsection{Automaton Evolution}

Let $G = (V,E)$ denote the square lattice graph, where $V$ is the set of lattice points 
and $E$ is determined by the chosen neighborhood structure. 
For $p \in V$, we denote by 
\[
Ne(p) = \{ g \in V : (p,g) \in E \}
\]
the set of neighboring vertices of $p$.

Consider a sequence
\[
\mathcal{S} =(k_1, k_2, k_3,\ldots), 
\qquad k_i \in \{2,3,4,\ldots\}.
\]

The discrete graph Laplacian applied to $\nu_i$ is defined by
\begin{equation}
	(\Delta \nu_{i})(p) = \sum_{g \in Ne(p)} \big( \nu_i(g) - \nu_i(p) \big)
	\label{eq:Laplacian}
\end{equation}
defined over $\mathbb{Z}$ prior to modular reduction.
The evolution rule is given by
\[
\nu_{i+1} = (\Delta \nu_i) \mod k_i,
\]
where $k_i$ is the $i$-th element of $\mathcal{S}$.

Here:
\begin{itemize}
	\item $k_i$ is the modulus (binary for $k_i=2$, ternary for $k_i=3$, etc.),
	\item $i$ denotes the discrete time (iteration) step,
	\item $Ne(p)$ defines the neighborhood structure (e.g., von Neumann or diagonal).
\end{itemize}

Thus, the automaton evolves by repeated application of the graph Laplacian, 
followed by reduction modulo $k_i$, for $i=1,2,3,\ldots$. 

\begin{figure}[h!]
	\centering
	\includegraphics[width=0.56\textwidth]{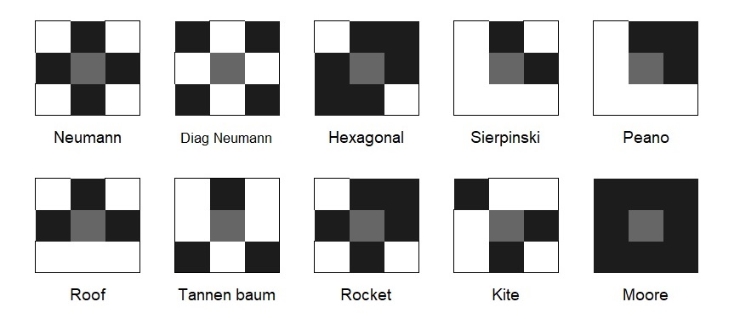}
	\caption{Examples of neighborhoods.}
	\label{fig:neighborhoods}
\end{figure}

\subsection{Modular Arithmetic and Iterative Extensions}

The introduction of the sequence $\mathcal{S}$ allows us to extend the construction 
beyond binary (black-and-white) figures, which correspond to the case when 
$\mathcal{S}$ is constant, i.e. $\mathcal{S} = (2,2,2,\ldots)$.

In particular, for constant sequences $\mathcal{S} = (n,n,n,\ldots)$ we obtain:

\begin{itemize}
	\item Ternary iteration ($n=3$): a three-state structure with increased combinatorial complexity.
	\item Quaternary and higher ($n>3$): multi-state structures with $n$ distinct states.
\end{itemize}

In this paper, we focus on sequences $\mathcal{S}$ of the form
\begin{equation}
	k_i =
	\begin{cases}
		n, & \text{if } i \equiv 2 \pmod{4}, \\
		2, & \text{otherwise},
	\end{cases}
\end{equation}
and the corresponding family of figures generated by the Laplacian dynamics 
will be referred to as \emph{$2n22$-style figures}\footnote{Other families of construction sequences are studied in 	\cite{NowakKepczyk2025}.}

In particular:
\begin{itemize}
	\item 2222-style: purely binary evolution ($n=2$),
	\item 2322-style: ternary interference at every fourth step ($n=3$).
\end{itemize}

\begin{figure}[h!]
	\centering
	\includegraphics[width=0.5\textwidth]{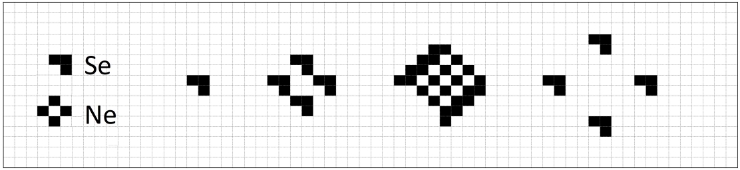}
\caption{Binary ($2222$-style) evolution from a three-point seed under the von Neumann neighborhood.  The upper left panel shows the initial seed and the assumed neighborhood structure.  The subsequent panels display the $0$-th (initial configuration), first, second, and third iterations  generated by the modular Laplacian dynamics with $k_i=2$, illustrating the initial expansion induced by the discrete Laplacian.\label{Fig3}}
\end{figure}

\section{Binary style figures and their properties}

 Prior work on binary style figures by Suzuki and Mae\-ga\-i\-to uncovered visually striking configurations, including snowflake-like, butterfly-like, and Persian carpet-like motifs \cite{Aiba}. Theoretical developments, including fixed-point theorems, periodicity results, and connections to binomial and trinomial sequences have been found \cite{Hadlich}. 
 
 We shall investigate binary style figures via many seeds and neighbourhoods. First we shall describe characteristic feature of binary figures and the suprising appearance of Sierpinski-like triangles.

To analyze connectedness, we shall say the figure is $k$-steps away from connectedness, $k> 0$, $k\in \mathbb{Z}$, if, in order to make it connected, it is required to add paths between its connected components, the longest being of length $k$. 

A connected figure ($A$), figures: 1-step away (figure $B$),  and 2-steps away from connectedness (figure $C$)  are shown below.\bigskip

\begin{center}
	\includegraphics[width=0.38\textwidth]{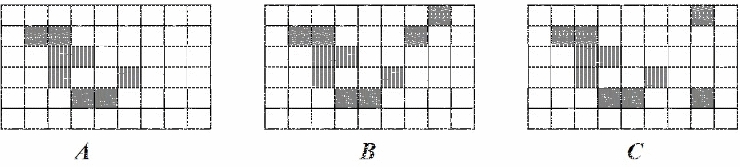}
\end{center}\smallskip


\noindent{\bf Proposition.}
{\it
	Let the binary Laplacian dynamics be defined with the Diag-Neumann neighborhood 
	and let the initial configuration be finitely supported. 
	Then for every $k\ge 1$, the configuration at iteration $i=8k$ 
	consists of disjoint translated copies of the initial seed. 
	Moreover, these copies are separated by a positive gap 
	(which equals $13$ in the $3\times3$ case considered here).
}
\bigskip

\noindent An illustration of this proposition is shown in Fig. \ref{Fig3}. The sketch of the proof will be given in the Appendix. 

\paragraph{Remark on other neighborhoods.}
The above argument is carried out for the Diag-Neumann neighborhood. 
Numerical experiments indicate that similar replication phenomena 
occur for other neighborhood structures as well. 
A general operator-based proof for arbitrary neighborhoods 
is beyond the scope of the present paper.

\subsection{Density Measure} 
The density of a figure provides a quantitative measure of its occupancy on the lattice and offers insights into its connectivity properties. A low density often correlates with fragmentation, where the figure consists of multiple disconnected components.

\medskip

\noindent
{\bf Definition.} The density \(\rho(i)\) of a figure at iteration \( i \) is defined as the ratio of occupied cells to the total available lattice area:
\[
\rho(i) = \frac{\sum_{p \in \text{lattice}} \text{sgn}(\nu(p))}{(3 + 2i)^2}.
\]

\noindent
Since binary figures are observed to become a spread of seeds at iterations \( i=8k, \, k=1,2,\dots \), their density exhibits sharp reductions at these steps (see Fig. \ref{graph1} (a)). This pattern is consistent with our analytical findings, suggesting that binary figures experience periodic dissociation into separated copies of their initial seeds.

At iteration \( i=8 \), density does not exceed:
\[
\rho(8) \leq \frac{36}{19^2} \approx 0.1.
\]

\subsection{Sierpinski-like Triangles} 
While binary constructions at iterations \( 8k \), \( k=1,2,\ldots \) reduce the figure to a spread of seeds, a distinct behavior emerges in a specific subsequence of iterations:
\[
i = 2^{k+1} - 1, \quad k=0,1,2,\ldots.
\]
At these steps, the figures reach their maximum local density (see Fig. \ref{graph1}) and, remarkably, form Sierpinski-like triangles—a behavior observed independently of the initial seed configuration.

An example of Sierpinski-like triangle growth is illustrated in Fig. \ref{SierpinskiLike.jpg}(a). The estimated fractal dimension of the figure at the sixth iteration, obtained via the box-counting method, is approximately:
\[
D_f \approx 1.51.
\]

\begin{figure}[h!]\begin{center}
		\includegraphics[width=0.9\textwidth]{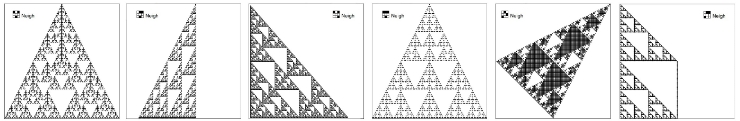}\end{center}	
	\caption{Sierpinski-like triangle constructions obtained as 2222-style figures at iterations \( i = 2^{k+1}-1 \), where \( k=0,1,2,\ldots \). The same neighborhood applied to different seeds produces structurally similar results, though the fractal dimension (via box-counting) varies depending on the seed.}
	\vspace{1em} 
	\label{SierpinskiLikeColl.jpg}
\end{figure}

\begin{figure}[h!]\begin{center}
		(a)	\includegraphics[width=0.9\textwidth]{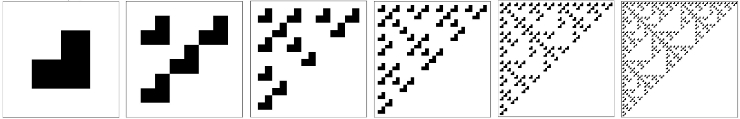}\\
		(b) \includegraphics[width=0.9\textwidth]{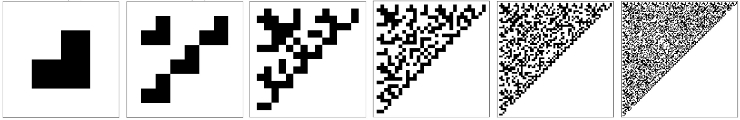}
	\end{center}	
	\caption{Comparison of 2222-style and 2322-style constructions from the same seeds and neighborhoods at the iteration sequence \( i = 2^{k+1}-1 \), where \( k=0,1,\ldots ,5\). The 2322-style construction exhibits more irregularity and disorder compared to the structured Sierpinski-like growth of 2222-style figures. }
	\vspace{1em} 
	\label{SierpinskiLike.jpg}
\end{figure}

A broader classification of five types of Sierpinski-like figures, each arising from different seeds but following the same iterative process, is shown in Fig. \ref{SierpinskiLikeColl.jpg}.

Additionally, we compare the 2222-style and 2322-style figures generated at the same subsequence of iterations in Fig. \ref{SierpinskiLike.jpg}. Notably, while 2222-style figures form clear, structured Sierpinski-like patterns, the 2322-style construction introduces greater structural complexity, leading to an apparently more chaotic formation.

\section{Unique Properties of $2322$-Style Figures}

The 2322-style construction, unlike purely binary (2222) figures, introduces ternary interactions, resulting in distinctive structural and dynamical properties. These include:
\begin{itemize}
	\item Lower density variance across iterations.
	\item Minimal connectivity loss (at most two steps away from full connectivity).
	\item Non-repetitive structural variation, distinct from the periodicity seen in 2222 figures.
\end{itemize}

\begin{figure}[h!]\begin{center}
		(a)\includegraphics[width=0.54\textwidth]{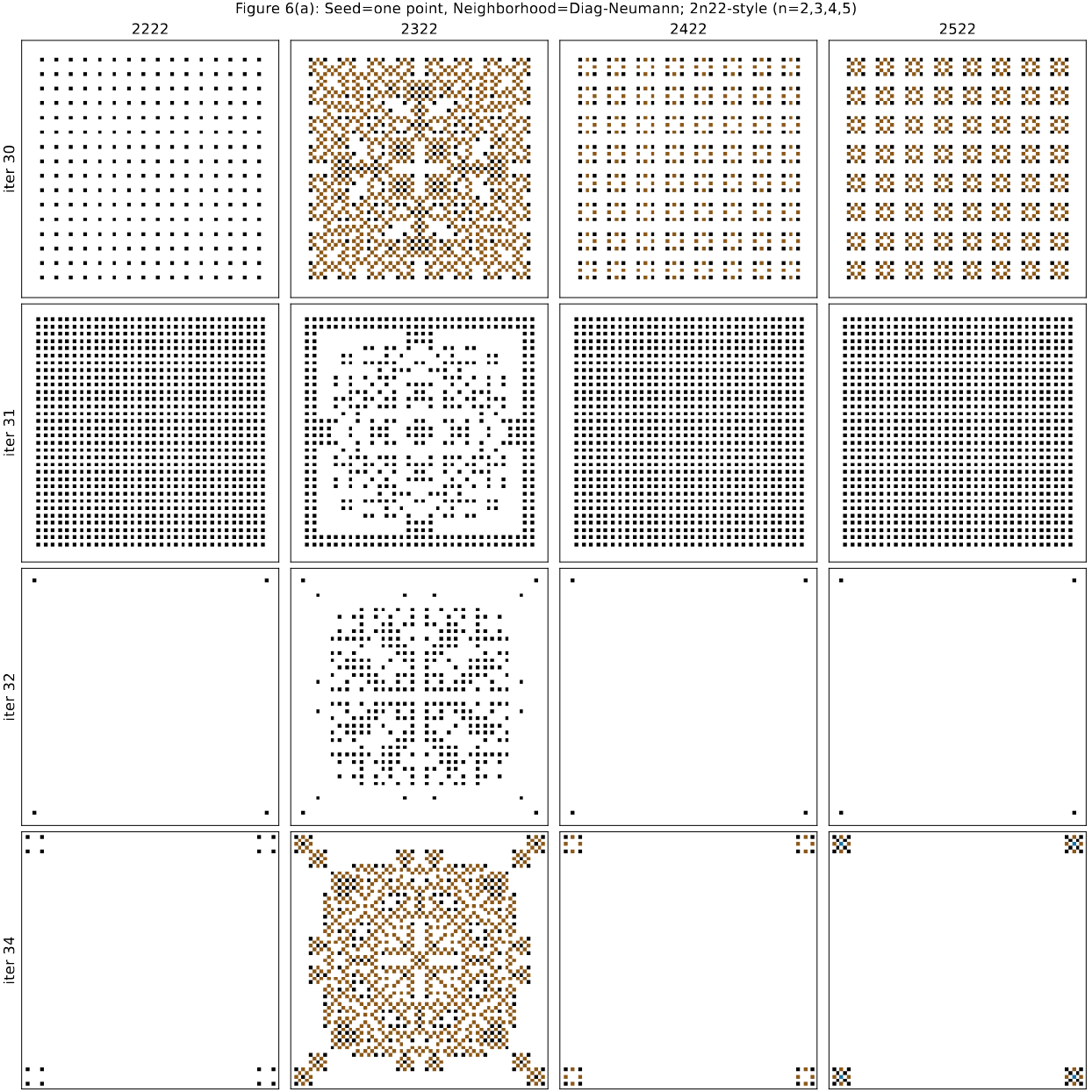}\\
		(b)\includegraphics[width=0.54\textwidth]{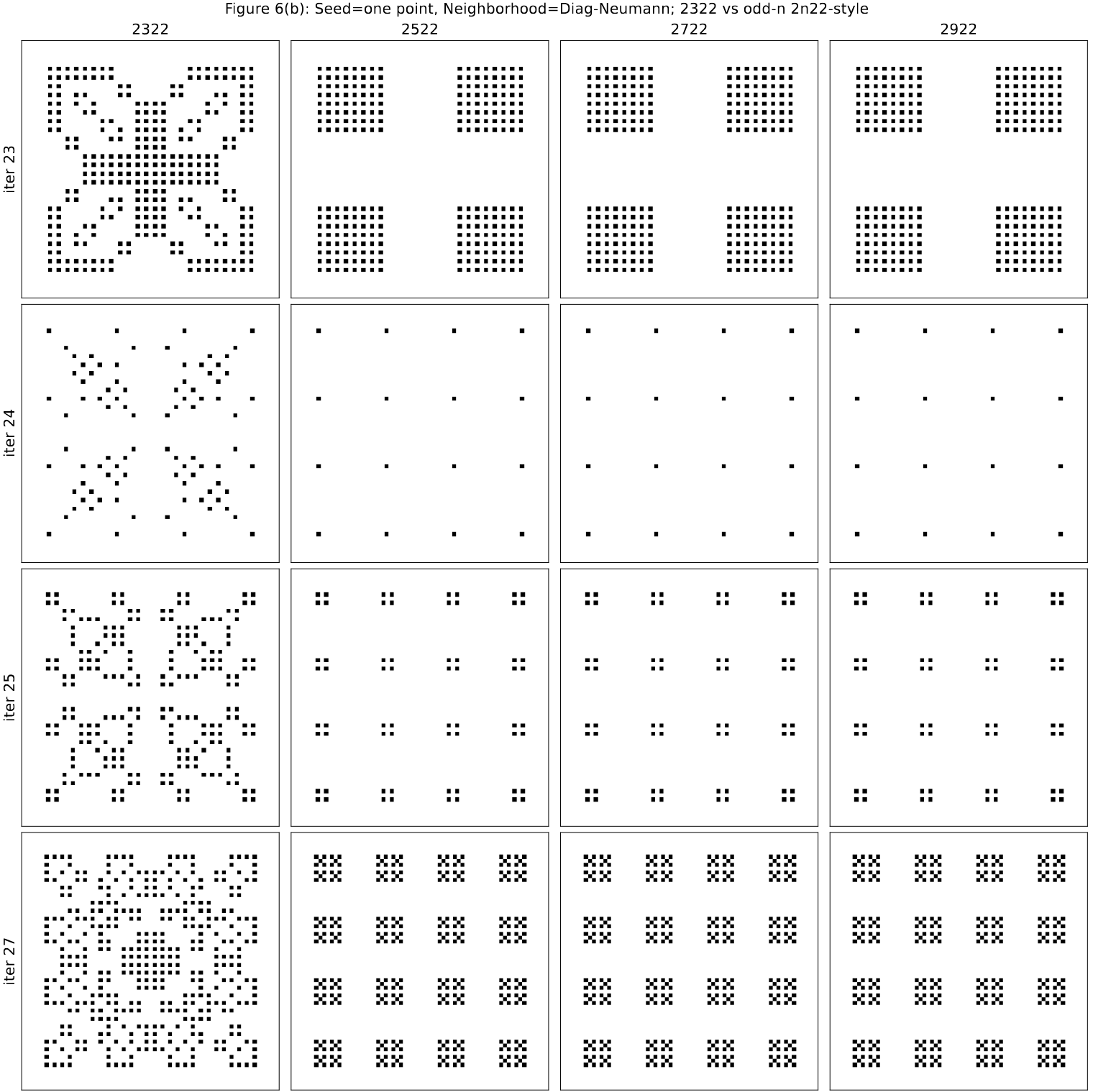}
	\end{center}
\caption{Seed: single point; diagonal Neumann neighborhood. \\
	(a) Comparison of $2n22$ figures for $n=2,3,4,5$ (left to right; iterations 30–32,34). 
	The 2322 case (second column) is structurally distinct. 
	Fractal dimensions for 2222 (resp. 2322) are 1.31, 1.66, 0.91, 1.14 
	(resp. 1.78, 1.62, 1.43, 1.62). \\
	(b) Comparison of 2322 with $2n22$ for odd $n=3,5,7,9$ (iterations 23–25,27), 
	showing persistent structural distinctness.}
	\label{ALL}

\end{figure}

At higher iterations, 2322-style figures retain complexity, while higher-order $2n22$-figures (for $n=5,7,9$) tend to exhibit more uniform growth patterns (Fig. \ref{ALL}).

\subsection{Connectivity and Structural Complexity}

\begin{figure}[h!]
	\centering
	\includegraphics[width=0.94\textwidth]{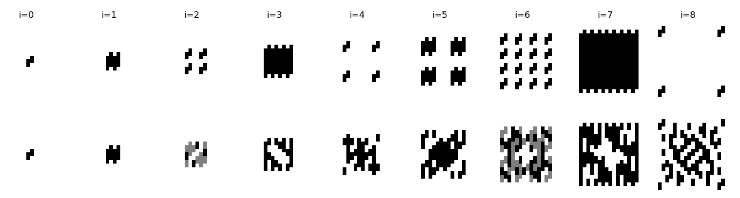}
	\caption{Comparison of 2222 and 2322 evolutions from the same seed. Unlike 2222 figures, 2322-style figures maintain structural connectivity, in contrast to the seed dissociation observed in 2222 at 8$i$-ths iterations.}
	\label{Prop2322}
\end{figure}

The distinct connectivity properties of 2322-style figures can be observed in their early evolution. Fig. \ref{Prop2322} illustrates how ternary interactions prevent the dissipation of seeds, which is characteristic of purely binary figures at iteration 8.

Among 2$n$22-style figures for $n$ odd, 2322 exhibits distinct structural behaviour:
\begin{itemize}
	\item It is the most connected, requiring at most two steps to fully connect all components.
	\item Figures for $n=5,7,9$ tend to lose structural diversity, converging into similar shapes differing mainly in color distribution.
	\item 2322 figures remain structurally distinct, retaining dynamic, non-repetitive features.
\end{itemize}

\subsection{Density Variability and Aperiodic Behavior}

\begin{figure}[h!]\begin{center}
		(a) \includegraphics[width=0.7\textwidth]{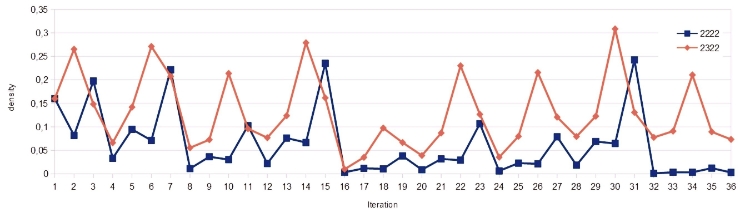} \\
		(b) \includegraphics[width=0.7\textwidth]{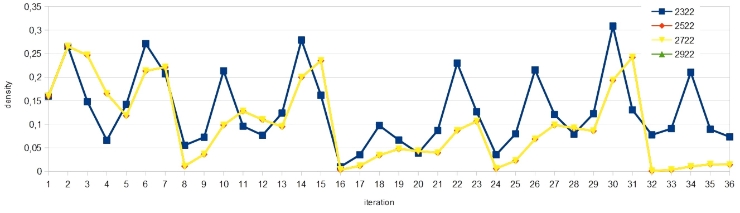} 
		\caption{Comparison of densities of figures over iterations, Seed: One Point, Neighborhood: Diag Neumann. \\(a) 2222-style and 2322-style figures over iterations.  
		Local minima at iterations $8k$ correspond to seed dissociation in 2222 figures, a behavior not observed in 2322 figures.\\
		(b) 2$n$22-style figures, $n=3,5,7,9$.
		Graphs for $n=5,7,9$ overlap, suggesting similar periodic behavior, whereas 2322 maintains a unique profile.\label{graph1}}
			\end{center}
\end{figure}

The density evolution of 2222 and 2322 figures over time is compared in Fig. \ref{graph1}. Notably, 2222 figures undergo periodic density reductions at iterations \( 8k \), $k=1,2,\dots$, corresponding to their spread-of-seeds phenomenon. This is absent in 2322-style figures, which exhibit a more consistent density profile.

Furthermore, we analyze density trends for odd-$n$ figures in Fig. \ref{graph1}(b). While figures for $n=5,7,9$ overlap, suggesting uniform periodicity, the 2322 figure remains structurally different.

\subsection{Connection to Aperiodic Sequences}

The non-periodic structural variation of 2322-style figures suggests connections to aperiodic mathematical structures, including:
\begin{itemize}
	\item Dekking’s Construction: A non-repetitive sequence model in combinatorial mathematics.
	\item Quasi-periodic tilings: Found in quasicrystals and aperiodic lattices.
	\item Fractal substitution sequences, seen in hierarchical dynamical systems.
\end{itemize}

These findings reinforce the idea that modular arithmetic sequences within discrete Laplacians can yield quasi-periodic, self-organizing patterns, with potential applications in mathematical modeling, materials science, and biological growth simulations.

\section{Statistical Comparison of 2222 and 2322 Sequences}

This section investigates the numerical properties of discrete Laplacian sequences generated by 2222-style and 2322-style constructions at a fixed lattice cell. The goal is to determine whether these structures exhibit long-range correlations, spectral properties, and fractal characteristics, providing potential links to Dekking’s Construction and other non-periodic sequences.

\subsection{Statistical Characteristics}

The entropy is computed in the Shannon sense,
\[
H = - \sum_{j} p_j \log p_j,
\]
where $p_j$ denotes the empirical frequency of the $j$-th state.
The fractal dimension is estimated using a standard box-counting method.

Table \ref{tab:numerical_stats} summarizes key statistical properties over 500 iterations. The entropy measures the degree of disorder, while variance and fractal dimension indicate complexity.

\begin{table}[h]
	\centering
	\begin{tabular}{|c|c|c|c|c|}
		\hline
		\textbf{Sequence Type} & \textbf{Entropy} & \textbf{Mean} & \textbf{Variance} & \textbf{Fractal Dimension} \\
		\hline
		2222-style & 0.1721 & 0.0275 & 0.0253 & 1.432 \\
		2322-style & 1.0517 & 0.3318 & 0.3425 & 1.867 \\
		\hline
	\end{tabular}
	\caption{Comparison of entropy, statistical measures, and fractal dimensions for 2222 and 2322 sequences over 500 iterations.}
	\label{tab:numerical_stats}
\end{table}

The higher entropy of the 2322 sequence indicates greater unpredictability compared to the periodic 2222 sequence. Similarly, the increased variance and higher fractal dimension suggest that 2322 structures exhibit greater structural complexity.

\subsection{Fourier Spectral Analysis}

\begin{figure}[h!]
	\centering
	\includegraphics[width=0.7\textwidth]{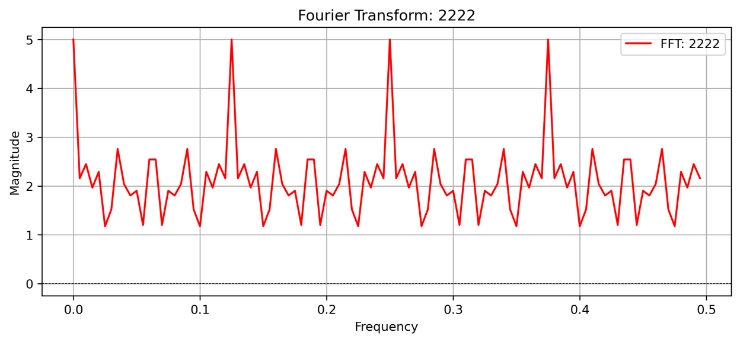}
	\caption{Fourier spectrum of the 2222 sequence. Sharp peaks indicate strong periodicity.}
	\label{fig:fourier_2222}
\end{figure}

\begin{figure}[h!]
	\centering
	\includegraphics[width=0.7\textwidth]{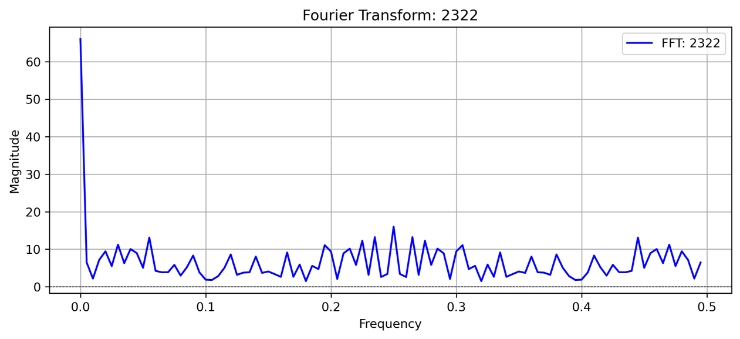}
	\caption{Fourier spectrum of the 2322 sequence. The broader spectrum suggests quasi-periodicity or hierarchical structure.}
	\label{fig:fourier_2322}
\end{figure}

Figures \ref{fig:fourier_2222} and \ref{fig:fourier_2322} present 
the amplitude spectra of the corresponding sequences.

The 2222 sequence exhibits pronounced discrete peaks 
located at regularly spaced frequencies, 
indicating strong periodic components 
and a well-defined fundamental frequency.

In contrast, the 2322 sequence does not display 
a single dominant periodic frequency. 
Instead, its spectrum consists of multiple distributed peaks 
with irregular spacing and varying amplitude, 
reflecting hierarchical spatial organization 
rather than strict periodic repetition.

Thus, while the 2222 case corresponds to a predominantly discrete spectrum, 
the 2322 sequence exhibits a more distributed spectral structure, 
consistent with non-periodic yet structured growth.

\subsection{Autocorrelation Analysis}

\begin{figure}[h!]
	\centering
	\includegraphics[width=0.7\textwidth]{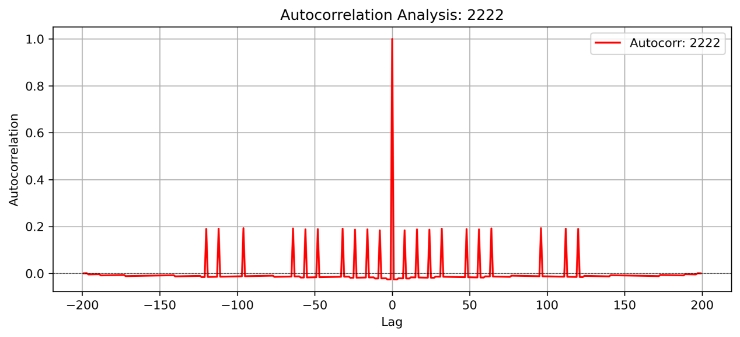}
	\caption{Autocorrelation function of the 2222 sequence. Periodic correlations indicate a structured, repeating pattern.}
	\label{fig:auto_2222}
\end{figure}

\begin{figure}[h!]
	\centering
	\includegraphics[width=0.7\textwidth]{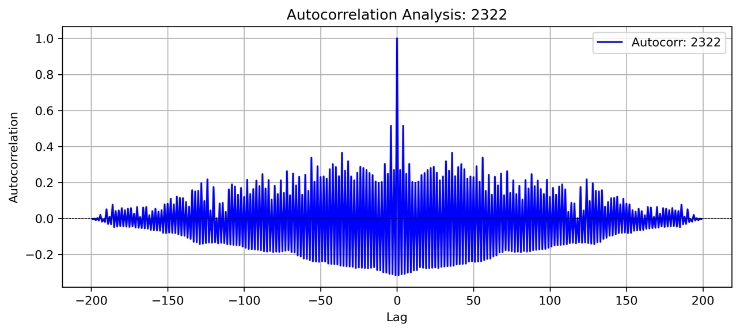}
	\caption{Autocorrelation function of the 2322 sequence. The irregular decay suggests non-trivial long-range correlations.}
	\label{fig:auto_2322}
\end{figure}

The autocorrelation function, shown in Figures \ref{fig:auto_2222} and \ref{fig:auto_2322}, provides further insight into sequence regularity. 

\begin{itemize}
	\item The 2222 sequence exhibits clear periodic correlations, confirming its structured, repeating nature—akin to the Sierpiński triangle.
	\item The 2322 sequence, in contrast, displays weaker long-range correlations and an irregular decay, suggesting a lack of strict periodicity.
\end{itemize}

This aperiodic structure suggests that 2322 figures belong to a distinct mathematical class, exhibiting properties found in quasicrystals, aperiodic tilings, and hierarchical dynamical systems.

\subsection{Connection to Dekking’s Construction}

Dekking’s Construction \cite{Dekking1979} generates non-periodic sequences that avoid Abelian repetitions. The 2322 sequence, with its high entropy, irregular autocorrelation, and absence of Fourier periodicity, shows structural similarities to this framework.

\begin{itemize}
	\item 2222 figures behave as classical fractals, with periodic behavior and predictable density fluctuations.
	\item 2322 figures exhibit aperiodic, self-organizing growth, leading to a higher fractal dimension and more complex scaling properties.
\end{itemize}

These findings suggest that discrete Laplacians could serve as a graphical realization of Dekking-like sequences, providing new insights into hierarchical self-similarity, quasi-periodicity, and aperiodic mathematical structures.

\subsection{Future Research Directions}

Future studies could explore:
\begin{itemize}
	\item Higher-dimensional extensions: Extending discrete Laplacian iterations to 3D lattices to study their higher-rank algebraic properties.
	\item Stochastic variations: Investigating randomized Laplacian growth to model natural fractal phenomena.
	\item Mathematical characterization: Developing a rigorous classification of ternary interference in modular growth patterns.
\end{itemize}

These results suggest that the 2322 sequence represents a novel class of aperiodic, fractal-like structures, meriting further investigation in both pure mathematics and applied physics.

\section{Conclusions and Perspectives}

We have studied modular Laplacian dynamics on two-dimensional square lattices,
with particular emphasis on non-constant sequences such as the $2322$ pattern.
In contrast to purely binary ($2222$) constructions,
which exhibit periodic behavior and seed dissociation at specific iterations,
the $2322$-style figures display reduced periodicity,
higher entropy, and increased fractal dimension.

These results indicate that introducing variable modular arithmetic
significantly modifies the long-term structural behavior
of discrete Laplacian iterations.

\subsection{Mathematical Implications}

\begin{itemize}
	\item \textbf{Non-periodic growth:}
	The structural complexity of $2322$-style figures
	suggests analogies with non-periodic substitution systems,
	including constructions related to Dekking-type sequences \cite{Dekking1979}.
	This correspondence is structural rather than formal.
	
	\item \textbf{Fractal characteristics:}
	Box-counting estimates indicate fractal dimensions
	in the approximate range $0.9$--$1.87$,
	consistent with the observed multi-scale self-similarity.
	
	\item \textbf{Higher-dimensional extensions:}
	Extending modular Laplacian dynamics
	to three-dimensional lattices and other graph structures
	remains a natural direction for further study.
\end{itemize}

\subsection{Physical and Computational Perspectives}

\begin{itemize}
	\item \textbf{Aperiodic order:}
	The reduced periodicity of $2322$-style figures
	shows structural similarities to systems exhibiting
	aperiodic order, such as quasicrystalline patterns \cite{GellMann1994}.
	This analogy concerns geometric organization
	rather than a direct physical model.
	
	\item \textbf{Fractal-inspired design:}
	Self-similar geometries are widely used in electromagnetic and sensor applications \cite{Antena,packages}.
	The controlled irregularity of the $2322$-style figures
	suggests potential relevance for scalable pattern generation.
	
	\item \textbf{Structured discrete models:}
	The non-uniform yet organized density patterns
	provide an abstract framework for studying
	hierarchical structures in graph-based systems.
	
	\item \textbf{Motivation from natural self-assembly:}
	Recent observations of Sierpiński-like protein structures 
	in biological systems \cite{Sendker,denovo,Singh}
	demonstrate that fractal growth patterns 
	may arise in natural self-assembly processes. 
	While the present work is purely mathematical, 
	such findings provide motivation for studying 
	discrete self-similar dynamics.
\end{itemize}

\subsection{Future Directions}

Possible directions for further investigation include:

\begin{itemize}
	\item Analytical characterization of non-binary modular sequences;
	\item Stochastic variants and perturbed modular iterations;
	\item Extensions to higher-dimensional lattices.
\end{itemize}

These directions aim to deepen the mathematical understanding
of modular Laplacian dynamics and their structural properties.

\appendix
\section*{Appendix: Proof Sketch}

We consider an initial seed contained in a $3 \times 3$ square and apply the Diag-Neumann neighborhood, defined by:

\[
Nei =
\begin{bmatrix}
	1 & 0 & 1 \\
	0 & 0 & 0 \\
	1 & 0 & 1
\end{bmatrix}
\]

Let the initial seed $F_0$ be:

\[
F_0 =
\begin{bmatrix}
	a & b & c \\
	d & \mathbf{e} & f \\
	g & h & i
\end{bmatrix}
\]

where $a, b, c, d, e, f, g, h, i \in \{0,1\}$ represent binary values.

We show that for $i=8k, k=1,2,3,\ldots$, the figure consists of a spread of its seeds and is at least 13-steps away from connectedness.

\paragraph{Step 1: First Iteration $F_1$}
Applying the Laplacian update:

\[
a^1_{x,y} = \sum_{(u,v) \in Ne(x,y)} \left( a^0_{x+u, y+v} - a^0_{x,y} \right) \mod 2
\]
yields:
\[
F_1 =
\begin{bmatrix}
	a & b & ac & b & c \\
	d & e & df & e & f \\
	ag & bh & \mathbf{acgi} & bh & ci \\
	d & e & df & e & f \\
	g & h & gi & h & i
\end{bmatrix}
\]
At this step, the pattern has expanded outward but remains connected.

\paragraph{Step 2: Second Iteration $F_2$}
Continuing with the same update rule, we obtain:

\[
F_2 =
\begin{bmatrix}
	F_0 & 0_{3\times1} & F_0 \\
	0_{1\times3} & 0 & 0_{1\times3} \\
	F_0 & 0_{3\times1} & F_0
\end{bmatrix}
\]

At this stage, the pattern doubles in size but already shows a structure of dispersed seeds.

\paragraph{Step 3: Third Iteration $F_3$}
Now we obtain:

\[
F_3 =
\begin{bmatrix}
	a & b & ac & b & ac & b & ac & b & c \\
	d & e & df & e & df & e & df & e & f \\
	ag & bh & acgi & bh & cdgi & bh & acgi & bh & ci \\
	d & e & df & e & df & e & df & e & f \\
	ag & bh & acgi & bh & cdgi & bh & acgi & bh & ci \\
	d & e & df & e & df & e & df & e & f \\
	g & h & gi & h & gi & h & gi & h & i
\end{bmatrix}
\]
This is the last iteration before modular 2 cancellation starts to take effect.

\paragraph{Step 4: Fourth Iteration $F_4$}
Applying the Laplacian update again, we obtain:

\[
F_4 =
\begin{bmatrix}
	F_0 & 0_{5\times3} & F_0 \\
	0_{5\times3} & 0_{5\times5} & 0_{5\times3} \\
	F_0 & 0_{5\times3} & F_0
\end{bmatrix}
\]

\paragraph{Steps 5-7: Expansion of Each $F_0$ Within $F_4$}
Each of the four $F_0$ components in $F_4$ will now expand similarly to $F_3$ at the seventh step. That means we are one step away from a fully expanded spread.

\paragraph{Step 8: Modular 2 Cancellation}

At iteration $i=8$, the only additional effect comes from overlaps 
between the expanding translated blocks. 
However, since the evolution is performed modulo $2$, 
any contributions appearing twice cancel out.

The translation offset at this stage equals $16$ in each coordinate direction. 
Since the initial seed is contained in a $3\times3$ square, 
its maximal diameter does not exceed $3$. 
Therefore, the minimal separation between the translated copies at iteration $8$ 
is given by
\[
16 - 3 = 13.
\]
Thus, the configuration $F_8$ consists of four disjoint translated copies 
of the initial seed, separated by a positive gap. 
This completes the sketch of the proof.
\qed

\section*{Declarations}

\begin{itemize}
	\item Funding: Not applicable.
	\item Conflicts of Interest: The author declares no conflicts of interest.
	\item Author Contributions: Ma\l gorzata Nowak-K\c epczyk conducted all research and analysis.
\end{itemize}

\bibliographystyle{unsrt}

\bibliography{references}
\end{document}